\documentclass[sn-mathphys-num]{sn-jnl}

\usepackage{amssymb}
\usepackage{graphicx}%
\usepackage{multirow}%
\usepackage{amsmath,amssymb,amsfonts}%
\usepackage{amsthm}%
\usepackage{mathrsfs}%
\usepackage[title]{appendix}%
\usepackage{xcolor}%
\usepackage{textcomp}%
\usepackage{manyfoot}%
\usepackage{booktabs}%
\usepackage{algorithm}%
\usepackage{algorithmicx}%
\usepackage{algpseudocode}%
\usepackage{listings}%

\theoremstyle{thmstyleone}%
\newtheorem{theorem}{Theorem}
\newtheorem{lemma}{Lemma}

\theoremstyle{thmstyletwo}%
\newtheorem{remark}{Remark}%

\newtheorem{problem}{Problem}%

\theoremstyle{thmstylethree}%
\newtheorem{definition}{Definition}%

\newcounter{alphabet}

\newenvironment{Thm}[1][]{\refstepcounter{alphabet}%
	\bigskip%
	\noindent%
	{\bf Theorem \Alph{alphabet}}%
	\ifthenelse{\equal{#1}{}}{}{ (#1)}%
	{\bf .} \itshape}{\vskip 8pt}

\newcommand{\IN}{{\mathbb N}}
\newcommand{\IC}{{\mathbb C}}
\newcommand{\ID}{{\mathbb D}}

\def\be{\begin{equation}}
\def\ee{\end{equation}}
\newcommand{\blem}{\begin{lemma}}
\newcommand{\elem}{\end{lemma}}
\newcommand{\bthm}{\begin{theorem}}
	\newcommand{\ethm}{\end{theorem}}
\newcommand{\bcor}{\begin{cor}}
	\newcommand{\ecor}{\end{cor}}
\newcommand{\beg}{\begin{exam}}
	\newcommand{\eeg}{\end{exam}}
\newcommand{\begs}{\begin{examples}}
	\newcommand{\eegs}{\end{examples}}
\newcommand{\bdefe}{\begin{definition}}
	\newcommand{\edefe}{\end{definition}}
\newcommand{\bprob}{\begin{problem}}
	\newcommand{\eprob}{\end{problem}}
\newcommand{\bques}{\begin{ques}}
	\newcommand{\eques}{\end{ques}}
\newcommand{\bei}{\begin{itemize}}
	\newcommand{\eei}{\end{itemize}}
\newcommand{\bcon}{\begin{conj}}
	\newcommand{\econ}{\end{conj}}
\newcommand{\bop}{\begin{op}}
	\newcommand{\eop}{\end{op}}

\newcommand{\bas}{\begin{assertion}}
	\newcommand{\eas}{\end{assertion}}

\newcommand{\bfa}{\begin{fact}}
	\newcommand{\efa}{\end{fact}}

\newcommand{\bca}{\begin{ca}}
	\newcommand{\eca}{\end{ca}}

\newcommand{\bst}{\begin{step}}
	\newcommand{\est}{\end{step}}

\newcommand{\bsca}{\begin{sca}}
	\newcommand{\esca}{\end{sca}}

\newcommand{\bcl}{\begin{cl}}
	\newcommand{\ecl}{\end{cl}}

\newcommand{\bmlem}{\begin{mlem}}
	\newcommand{\emlem}{\end{mlem}}

\newcommand{\bscl}{\begin{scl}}
	\newcommand{\escl}{\end{scl}}

\newcommand{\bcons}{\begin{conjs}}
	\newcommand{\econs}{\end{conjs}}

\newcommand{\bprop}{\begin{prop}}
	\newcommand{\eprop}{\end{prop}}

\newcommand{\br}{\begin{remark}}
	\newcommand{\er}{\end{remark}}
\newcommand{\brs}{\begin{remarks}}
	\newcommand{\ers}{\end{remarkss}}
\newcommand{\bo}{\begin{obser}}
	\newcommand{\eo}{\end{obser}}
\newcommand{\bos}{\begin{obsers}}
	\newcommand{\eos}{\end{obsers}}
\newcommand{\bpf}{\begin{proof}}
	\newcommand{\epf}{\end{proof}}
\newcommand{\ba}{\begin{array}}
	\newcommand{\ea}{\end{array}}
\newcommand{\beq}{\begin{eqnarray}}
	\newcommand{\beqq}{\begin{eqnarray*}}
		\newcommand{\eeq}{\end{eqnarray}}
	\newcommand{\eeqq}{\end{eqnarray*}}

\newcommand{\ds}{\displaystyle}

\begin{document}

\title 
[Multidimensional Bohr inequalities involving Schwarz mappings]{Multidimensional analogues of the refined versions of Bohr inequalities involving Schwarz mappings}


\author[1]{\fnm{Shanshan} \sur{Jia}}\email{2071831393@qq.com}
\equalcont{All authors contributed equally to this work.}
\author*[1]{\fnm{Ming-Sheng} \sur{Liu}}\email{liumsh65@163.com}
\equalcont{All authors contributed equally to this work.}

\author[2,3]{\fnm{Saminathan} \sur{Ponnusamy}}\email{samy@iitm.ac.in}
\equalcont{All authors contributed equally to this work.}

\affil*[1]{\orgdiv{School of Mathematical Sciences}, \orgname{South China Normal University}, \orgaddress{\street{Guangzhou}, \city{Guangdong}, \postcode{510631}, 
\country{China}}}


\affil[2]{\orgdiv{Department of Mathematics}, \orgname{Indian Institute of Technology Madras},
\orgaddress{
\city{Chennai}, \postcode{600 036}, \state{Tamilnadu}, \country{India}}}

\affil[3]{\orgdiv{Moscow Center of Fundamental and Applied Mathematics},
\orgname{Lomonosov Moscow State University},
\orgaddress{
\city{Moscow},
\country{Russia}}
}


\abstract{Our first aim of this article is to establish several new versions of refined Bohr inequalities for bounded analytic functions in the unit disk involving Schwarz functions. Secondly, 
  we obtain several new multidimensional analogues of the refined Bohr inequalities for bounded holomorphic mappings on the unit ball in a complex Banach space involving higher dimensional Schwarz mappings. All the results are proved to be sharp.}




\maketitle

\pagestyle{myheadings}	
\markboth{S.S. Jia, M.S. Liu and S. Ponnusamy}{Multidimensional Bohr inequalities involving Schwarz mappings}
	
\section{Introduction and Preliminaries}\label{HLP-sec1}
	

Let $H_{\infty}$ denote the class of all bounded analytic functions $f$ in the unit disk $\ID= \{z\in \IC:\,|z|<1\}$ equipped with the topology of uniform convergence on compact subsets of $\ID$ with the supremum norm $||f||_{\infty}:=\mathop{\sup}\limits_{z\in \ID}|f(z)|$.

\subsection{Classical Bohr inequality and its recent implications}
A remarkable discovery of Herald Bohr \cite{HB1914} in 1914 states that if  $f\in H_{\infty}$, then 
$$
	B_0(f,r):= |a_0|+\sum_{n=1}^{\infty} |a_n| r^n\leq ||f||_{\infty} ~\mbox{ for $0\leq r\leq 1/6$,}
$$
	where $a_n=f^{(n)}(0)/n!$ for $n\geq 0$. But subsequently later, Riesz, Shur and Wiener, independently proved
that this inequality holds for the wider range $0\leq r\leq 1/3$, and the number $1/3$ cannot be improved.
This inequality is known as the classical Bohr inequality and the constant $1/3$ is the Bohr radius
as the function $\varphi_a(z)=(a-z)/(1-\overline{a} z)$ ($|a|<1$) from the group of automorphisms of the unit disk
shows. The paper of Bohr \cite{HB1914} indeed contains the proof of Wiener showing that the Bohr radius is $1/3$. Several other proofs of it are available in the literature. See the survey articles \cite{YAM2017,SRG2018} and the references therein for the plane case.
It is worth pointing out that there is no extremal function in $H_\infty$ such that the Bohr radius is precisely $1/3$ (cf. \cite{SAA2019}, \cite[Corollary 8.26]{SRG2018} and \cite{IRK2017}).
	Several extensions, modifications and improvements of Bohr's theorem in different settings may be obtained from many recent papers. See \cite{YAM2013, SAA2019,VA2022, CB2004, BhowDas-18, BomBor-04,SKS2022} and the references therein. For example, some other improved versions of Bohr's inequality may be obtained from \cite{IRK2018} and Huang et al. \cite{HY2020} established new sharp versions of Bohr-type inequalities for functions from  $f\in H_{\infty}$ by allowing Schwarz function in place of the initial coefficients in the power series representations of $f$ about origin. In the recent years, many scholars have extended this problem to functions of several complex variables (cf. \cite{LA2000,SK2023,SKM2024,RY2023} and the recent monograph by Defant et al. \cite{DGMP19}). In 1997, Boas and Khavinson \cite{HB1997} defined $n$-dimensional Bohr radius for the family of holomorphic functions bounded by $1$ on the unit polydisk, which inspired a lot in the recent years, including Bohr's phenomena for multidimensional power series. Djakov and Ramanujan \cite{PBD2000} continued the investigation further, including the discussion on $p$-Bohr radius (see also \cite{IRK2019}, where the conjecture of Djakov and Ramanujan has been solved affirmatively). In 2021, Liu and Ponnusamy \cite{MS2021} have established several versions of multidimensional analogues of improved Bohr's inequality.

One of our aims in this paper is to generalize or to improve multidimensional analogues of the refined Bohr inequalities involving Schwarz functions
and to prove that these refinements are  best possible. 

\subsection{Basic Notations}\label{HLP-sec1-1}
Let ${\mathcal B} = \{f\in H_\infty :\, \|f\|_\infty \leq 1 \}$ and, for $m\in \mathbb{N}:=\{1,2,\ldots \}$, let
	\beqq
	{\mathcal B}_m=\{\omega \in {\mathcal B}:\, \omega (0)= \cdots =\omega ^{(m-1)}(0)=0 ~\mbox{ and }~ \omega ^{(m)}(0)\neq 0 \} 
	\eeqq
so that ${\mathcal B}_1=\{\omega \in {\mathcal B}:\, \omega (0)=0 ~\mbox{ and }~ \omega '(0)\neq 0 \}.$
	Also, for $f(z)=\sum_{n=0}^{\infty} a_{n} z^{n}\in {\mathcal B}$ and $f_0(z):=f(z)-f(0)$, we let (as in \cite{PVW2021})
	$$B_{N}(f,r) :=  \sum_{n=N}^{\infty} |a_n| r^n ~~\mbox{for $N\geq 0$,} ~\mbox{ and }~
	\|f_0\|_r^2  :=   \sum_{n=1}^{\infty}\left|a_{n}\right|^{2} r^{2 n}~,
	$$
	and in what follows we introduce
	$$
A(f_0,r):=\left (\frac{1}{1+|a_0|}+\frac{r}{1-r}\right )\|f_0\|_r^2,
	$$
	which  helps to reformulate refined classical Bohr inequalities.

	\subsection{Recently known refined Bohr-type inequalities}
	Recently, Chen  et al. \cite{CKX2023}, Huang  et al. \cite{HY2020} and Liu  et al. \cite{LG2021} proved several Bohr-type inequalities. We first  recall a few of these results.

	\begin{Thm}\label{Theo-A}\cite[Theorem 2]{HY2020}
		Suppose that $f(z)=\sum_{n=0}^{\infty} a_{n} z^{n} \in\mathcal{B}$,  and $\omega \in\mathcal{B}_m$ for some $m \in \IN$.
		Then we have
		\begin{eqnarray*}
			B_0(f,r)+ A(f_0,r) + \left|f\left(\omega(z)\right)-a_{0}\right| \leq 1
		\end{eqnarray*}
		for $r \in\left[0, \zeta_{m}\right] $, where $\zeta_{m}$ is the unique root in $(0, 1/3]$ of the equation
		\begin{equation*}
			r^{m}(3-5r)+3 r -1=0,
		\end{equation*}
		or equivalently, $3r^m +2\sum_{k=1}^{m}r^{k}-1=0.$ The upper bound $\zeta_{m}$ cannot be improved.
		
		Moreover,
		\begin{eqnarray*}
			|a_{0}|^2+B_1(f,r)+ A(f_0,r)+ |f(\omega(z))-a_{0}| \leq 1
		\end{eqnarray*}
		for $r \in\left[0, \eta_{m}\right] $, where $\eta_{m}$ is the unique root in $(0, 1/2]$ of the equation
		\begin{equation*}
			r^{m}(2-3r) +2 r-1=0, 
		\end{equation*}\
		or equivalently, $2r^m +\sum_{k=1}^{m}r^{k}-1=0.$
		The upper bound $\eta_{m}$ cannot be improved.
	\end{Thm}
	
Several special cases of it are remarked in \cite[Remark 2]{HY2020}.

	\begin{Thm}\label{Theo-B}\cite[Theorem 6]{CKX2023}
		Suppose that $f(z)=\sum_{n=0}^{\infty} a_{n} z^{n} \in\mathcal{B}$, $\omega \in\mathcal{B}_m $, $p\in(0,2],\, m,\, q\in\IN,\, q\geq 2$, and $0<m<q$. Then for arbitrary $\lambda\in (0,\infty)$, we have
		\begin{equation*}
		    |f(\omega(z))|^p+\lambda\sum_{k=1}^{\infty}|a_{qk+m}|r^{qk+m}\leq 1\quad \mbox{for}\quad r\leq R^{p}_{\lambda,q,m},
		\end{equation*}
		where $R^{p}_{\lambda,q,m}$ is the minimal positive root of the equation $\Psi(r)=0$ in the interval $[0,1]$, where
		$$\Psi(r):=2\lambda\frac{r^{q+m}}{1-r^q}-p\frac{1-r^m}{1+r^m}.
		$$
		In the case when $\Psi(r)>0$ in some interval $(R^{p}_{\lambda,q,m}, R^{p}_{\lambda,q,m}+\varepsilon)$, the number $R^{p}_{\lambda,q,m}$ cannot be improved.
	\end{Thm}

	\begin{Thm}\label{Theo-C}\cite[Theorem 5]{CKX2023}
	    Suppose that $f(z)=\sum_{n=0}^{\infty} a_{n} z^{n} \in\mathcal{B}$,  $\omega \in\mathcal{B}_m$ for some $m \in \IN$ and $p\in(0,2]$. Then for arbitrary $\lambda\in (0,\infty)$, we have
		\begin{eqnarray*}
			\left|f\left(\omega(z)\right)\right|^p+\lambda[B_1(f,r)+ A(f_0,r)]\leq 1
		\end{eqnarray*}
		for $0\leq r\leq R_{\lambda,m}$, where $R_{\lambda,m}$ is the best possible and it is the minimal positive root in $[0, 1]$ of the equation
		\begin{equation*}
			2\lambda\frac{r}{1-r}-p\frac{1-r^m}{1+r^m}=0.
		\end{equation*}
	\end{Thm}

\begin{Thm}\label{Theo-D}\cite[Corollary 1]{LG2021}
If $f\in {\mathcal B}$ and $f(z)=\sum_{n=0}^{\infty} a_{pn} z^{pn}$ for some  $p\in\mathbb{N}$, then the following sharp inequalities hold:
\begin{enumerate}
	\item[{\rm (a)}] $\ds \sum_{n=0}^{\infty}|a_{pn}|r^{pn}+ \left(\frac{1}{1+|a_0|}
+\frac{r^p}{1-r^p}\right)\sum_{n=1}^{\infty} |a_{pn}|^2 r^{2pn}\leq1$ \quad \text{for}~~$\ds r\leq\frac{1}{\sqrt[p]{2+|a_0|}}.$
	\item[{\rm (b)}] $\ds |a_0|^2+\sum_{n=1}^{\infty}|a_{pn}|r^{pn}+ \left(\frac{1}{1+|a_0|}
+\frac{r^p}{1-r^p}\right)\sum_{n=1}^{\infty} |a_{pn}|^2 r^{2pn}\leq1$ \quad \text{for}~~$\ds r\leq\frac{1}{\sqrt[p]{2}}.$
\end{enumerate}
\end{Thm}

 \begin{Thm}\label{Theo-E}\cite[Theorem 5]{HY2020}
 Suppose that $f(z)=\sum_{n=0}^{\infty} a_{n} z^{n} \in\mathcal{B}$,  $a:=|a_0|$ and $\omega \in\mathcal{B}_m$ for some $m\geq 1$.  We have the following:
\begin{enumerate}
\item[{\rm (1)}] If $m=1$, then we have
\be
 I_f(z):=B_0(f,r)+\frac{1+ar}{(1+a)(1-r)} \|f_0\|_r^2+\left|f\left(\omega(z)\right)-a_{0}\right|^{2} \leq 1 ~\mbox{ for $ |z|=r\leq 1/3 $}
\label{liu29}
\ee
if and only if $0\leq a \leq a^{*}=-5 +4\sqrt{2}\approx 0.656854$. The constant $1/3$ cannot be improved.

\item[{\rm (2)}] If $m\geq 2$, then \eqref{liu29} holds, and the constant $1/3$ cannot be improved.
\end{enumerate}
\end{Thm}

All the above results either generalize or refine many of the earlier known results as special cases in the planar case.

\subsection{New problems on multidimensional Bohr's inequality}
Let $X$ and $Y$ be complex Banach spaces with norms $||\cdot||_X$ and $||\cdot||_Y$, respectively. For simplicity, we omit the subscript for the norm when it is obvious from the context. Let $B_X$ and $B_Y$ be the open unit balls in $X$ and $Y$, respectively. If $X=\IC$, then $B_X=\ID=\{z\in \IC:\,|z|<1\}$ is the unit disk in $\IC$. For each $x\in X\backslash\{0\}$, we define
$$
   T(x)=\{T_x\in X^*:\,||T_x||=1,\,T_x(x)=||x||\},
$$
where $X^*$ is the dual space of $X$. Then the famous Hahn-Banach theorem implies that $T(x)$ is non empty.

\bdefe
Let $X$ and $Y$ be complex Banach spaces and $n\in \IN$. A mapping $P:\,X\rightarrow Y$ is called a homogeneous polynomial of degree $n$ if there exists a $n$-linear mapping $u$ from $X^k$ into $Y$ such that
$$P(x) = u(x,\ldots,x)$$
for every $x\in X$.


\edefe

Throughout our discussion, we denote the degree of a homogeneous polynomial by a subscript. We note that if $P_j$ is an $j$-homogeneous polynomial from $X$ into $Y$, there uniquely exists a symmetric $j$-linear mapping $u$ with $P_j(x) = u(x,\ldots,x)$.

A holomorphic function $f:\, B_X\rightarrow \ID$ can be expressed in a homogeneous polynomial expansion $f(z)=\sum _{j=0}^\infty P_j(z),z\in B_X$, where $P_j(z) = \sum _{|\alpha|=j}c_\alpha z^\alpha$ is a homogeneous polynomial of degree $j$, and $P_0(z) = f(0)$. 

A holomorphic mapping $\mu:\, B_X\rightarrow B_Y$ with $\mu(0) = 0$ is called a Schwarz mapping. We note that if $\mu_n$ is a Schwarz mapping such that $z = 0$ is a zero of order $n$ of $\mu_n$, $n \in \IN$, then the following estimation holds (see e.g. \cite[Lemma 6.1.28]{IG2003}):
\begin{eqnarray}\label{eq1.2}
||\mu_n(z)||_Y\leq ||z||^n_X,\quad z \in B_X.	
\end{eqnarray}

It is natural to raise the following.
	
\bprob\label{HLP-prob1}
Can we establish the improved versions of Theorems A-E, 
by involving two Schwarz functions?
\eprob
\bprob\label{HLP-prob2}
Can we establish the multidimensional versions of Theorems A-E, 
by involving two Schwarz mappings?
\eprob
	
Our goal in this article is to present affirmative answers to these two problems.
	
The paper is organized as follows. In Section \ref{HLP-sec2}, we present statements of our theorems and several remarks; in particular, Theorems \ref{NEW-lem4} to \ref{NEW-lem9} essentially provide an affirmative answer to Problem \ref{HLP-prob1}, and Theorems \ref{HLP-th1} to \ref{HLP-th6} provide an affirmative answer to Problem \ref{HLP-prob2}. In Section \ref{HLP-sec3}, we present some necessary lemmas which are needed for the proofs of our theorems. In Section \ref{HLP-sec4}, we supply the proofs of all the theorems.

	\section{Statement of Main Results, Remarks and Special Cases}\label{HLP-sec2}
	
We first generalize Theorem A 
as follows.

\bthm\label{NEW-lem4}
	Suppose that $f \in\mathcal{B}$  has the expansion $f(z)=\sum_{j=0}^{\infty} a_j z^j$ with $a=|f(0)| < 1$, $\omega_n \in\mathcal{B}_n$ for $n\geq 1$ and $R_1:=R_{k,m}^p$ is the unique root in $(0, 1)$ of the equation $\Psi_1(r)=0$, where
\begin{equation}
   	\Psi_1(r)=\frac{r^k}{1-r^k}+\frac{r^m}{1-r^m}-\frac{p}{2} \label{eq2.2}
\end{equation}
for some $k, m \in \IN$ and $p\in(0,2]$. Then we have
	\begin{eqnarray*}
		A_f(z) := a^p+B_1(f,|\omega_k(z)|) + 
A(f_0,|\omega_k(z)|)
+\left|f\left(\omega_m(z)\right)-a_0\right| \leq 1
	\end{eqnarray*}
for $|z|= r \leq  R_1$. 
	The number $R_1$ cannot be improved.
	\ethm 

\br

   We mention a few useful remarks and some special cases of Theorem \ref{NEW-lem4}.
\begin{enumerate}
\item We observed that \eqref{eq2.2} yields that the roots of $\Psi_1(r)=0$ are same when the corresponding sums $k+m$ of $k$ and $m$ are the same. That is, for each fixed $i,j \in \IN$, if $k_i+m_i=k_j+m_j$, then $R^p_{k_i,m_i}=R^p_{k_j,m_j}$.

\item If we set $k=1,\, \omega_k(z)=z$, $p=1$ and $p=2$ in Theorem \ref{NEW-lem4}, then we get Theorem A.

\item  If we set $p= k= m=1$ in \eqref{eq2.2}, then we get  $R^1_{1,1}=1/5$.
\item  If we set $p=k=1$ and $m=2$ in \eqref{eq2.2}, then $R^1_{1,2}=\frac{\sqrt{6}-1}{5} ~(\approx 0.289898<1/3)$ as the minimal root of the equation
    $$ -5r^{3}+3r^2+3 r -1=(1-r)(5r^{2}+2r -1)=0.
    $$
\item  If we set $p=2$ and $k= m=1$ in \eqref{eq2.2}, then we get  $R^2_{1,1}=1/3 $.
\item  If we set $k=1$ and $p=m=2$ in \eqref{eq2.2}, then we get $R^2_{1,2}=\frac{\sqrt{13}-1}{6} ~(\approx 0.434259 <1/2)$ as the minimal root of the equation
$$ (1-r)(3r^{2}+ r -1)=0.
$$
\item  If we allow $m\rightarrow \infty $  in \eqref{eq2.2}  (with $\omega_m(z)=z^m$ ), then $|f(\omega_m(z))|\rightarrow |f(0)|$ and $R_{k,\infty}^p=\sqrt[k]{\frac{p}{2+p}}$. In particular, $R^1_{1,\infty}=1/3$ and $R^2_{1,\infty}=1/2$.
\end{enumerate}
In Table \ref{tab1}, we include the values of $R^p_{k,m}$ for certain other values of $m\geq 3$, especially when
$p=1, 2$ and $k=1,3,5$.
\er

   \begin{table}[tbp]
\centering
\begin{tabular}{|l|l||l|l||l|l|}
\hline
$m$  &$R^1_{1,m}$  &$R^2_{1,m}$  &$R^1_{3,m}$  &$R^1_{5,m}$ \\
\hline
3  &0.318201  &0.469396  &0.584804  &0.647197 \\
\hline
4  &0.328083  &0.484925  &0.624100  &0.695544 \\
\hline
5  &0.331541  &0.492432  &0.647197  &0.724780 \\
\hline
10 &0.333326  &0.499757  &0.685896  &0.780637 \\
\hline
15 &0.333333  &0.499992  &0.692116  &0.795317 \\
\hline
20 &0.333333  &0.500000  &0.693159  &0.800196 \\
\hline
\end{tabular}
\vspace{8pt}
\caption{Values of $R_{k,m}^p$ for certain choices of $p,k,m$ in Theorem \ref{NEW-lem4}} \label{tab1}
\end{table}

Next, we present a generalization of Theorem B 
    involving another Schwarz function.

\bthm\label{NEW-lem5}
	Suppose that $f(z)=\sum_{j=0}^{\infty} a_j z^j\in\mathcal{B}$, $\omega_n \in\mathcal{B}_n$ for $n\geq 1$, $p\in(0,2],\, s,\, t\in\IN,\, s\geq 2$, and $0<t<s$ and $R_2:=R^{p,m,k}_{\lambda,s,t}$ is the unique root in the interval $(0,1)$ of the equation  $\Psi_2(r)=0$, where
   \begin{equation}
   	  \Psi_2(r)=2\lambda\frac{r^{k(s+t)}}{1-r^{ks}}-p\frac{1-r^m}{1+r^m}\label{eq2.4}
   \end{equation}
for some $k, m \in \IN$ and $\lambda\in (0,\infty)$. Then we have
\begin{equation*}
		B_f(z):= |f(\omega_m(z))|^p+\lambda\sum_{j=1}^{\infty}\left|a_{sj+t}\right|\left|\omega_k(z)\right|^{sj+t}\leq 1\quad \mbox{for}\quad |z|= r\leq R_2. 
\end{equation*}
	
In the case when $\Psi_2(r)>0$ in some interval $(R_2, R_2+\varepsilon)$, where $\Psi_2$ is given by \eqref{eq2.4}, the number $R_2$ cannot be improved.
\ethm	

	
If we set $\lambda=1$, $s=1$ and $t=0$ in \eqref{eq2.4}, then the roots of the equation  $\Psi_2(r)=0$ is same as the minimal positive root $\alpha_{k,m,p}$ in $(0,1)$ of the equation
    \begin{equation}
   	  \zeta_{k,m,p}(r):=2r^k(1+r^m)-p(1-r^m)(1-r^k)=0.\label{eq2.11}
    \end{equation}
   In Table \ref{tab2}, we list the  values of $\alpha_{k,m,p}$ for certain choices of $k$, $m$ and $p$.

    \br
   It is worth mentioning certain other special cases of Theorem \ref{NEW-lem5}.
    \begin{enumerate}

        \item  If we set $p=k=m=1$ in \eqref{eq2.11}, then it reduces to $r^2+4r-1=0$ which gives the root  $\alpha_{1,1,1}=\sqrt{5}-2\approx 0.236068$.
        \item  If we set $p=m=1$ and $k=2$ in \eqref{eq2.11}, then it reduces to
        $$0=r^3+3r^2+r-1=(r+1)(r^2+2r-1)
        $$
        which gives the root $\alpha_{2,1,1}=\sqrt{2}-1\approx 0.414214$.
	    \item  If we set $p=2$ and $k=m=1$ in \eqref{eq2.11}, then we get $\alpha_{1,1,2}=1/3$.
        \item  If we set $p=m=2$ and $k=1$ in \eqref{eq2.11},
            then we get $\alpha_{1,2,2}=\sqrt{2}-1\approx 0.414214$.
        \item  If we set $p=k=2$ and $m=1$ in \eqref{eq2.11}, then we get $\alpha_{2,1,2}=1/2$.
        \item  If we allow $m\rightarrow \infty $  in \eqref{eq2.11}  (with $\omega_m(z)=z^m$ ), then $|f(\omega_m(z))|\rightarrow |f(0)|$ and $\alpha_{k,\infty,p}=\sqrt[k]{\frac{p}{2+p}}$. In particular, $\alpha_{1,\infty,1}=1/3$ and $\alpha_{1,\infty,2}=1/2$.
        \item If we set $k=1$ and $\omega_k(z)=z$ in Theorem \ref{NEW-lem5}, then we get Theorem B. 
            \end{enumerate}

\er
\begin{table}[tbp]
\centering
\begin{tabular}{|l||l|l||l|l|l|}
\hline
$m$  &$\alpha_{1,m,1}$  &$\alpha_{1,m,2}$  &$\alpha_{5,m,1}$  &$k$ &$\alpha_{k,1,1}$ \\
\hline
1  &0.236068  &0.333333  &0.632447  &1  &0.236068 \\
\hline
2  &0.295598  &0.414214  &0.686395  &2  &0.414214 \\
\hline
3  &0.319053  &0.453398  &0.715894  &3  &0.516239 \\
\hline
4  &0.328197  &0.474627  &0.735303  &4  &0.583776 \\
\hline
5  &0.331555  &0.486389  &0.749217  &5  &0.632447 \\
\hline
10 &0.333326  &0.499516  &0.783683  &10 &0.759593 \\
\hline
15 &0.333333  &0.499985  &0.795743  &15 &0.816751 \\
\hline
20 &0.333333  &0.500000  &0.800252  &20 &0.850170 \\
\hline
\end{tabular}
\vspace{8pt}
\caption{Values of $\alpha_{k,m,p}$ for certain choices of $k,m, p$ in Theorem \ref{NEW-lem5}} \label{tab2}
\end{table}

	Now, we state a new version of  Theorem C 
involving Schwarz functions.

\bthm\label{NEW-lem6}
	Suppose that $f(z)=\sum_{j=0}^{\infty} a_j z^{j} \in\mathcal{B}$,  $\omega_n \in\mathcal{B}_n$ for $n\geq 1$, and $p\in(0,2]$, $R_3: =R_{\lambda,p}^{k,m}$ is the unique root in $(0, 1)$ of the equation $\Psi_3(r)=0$, where
   \begin{equation}
   	\Psi_3(r)=p\frac{1-r^m}{1+r^m}-2\lambda\frac{r^k}{1-r^k}\label{eq2.6}
   \end{equation}
   for some $m, k\in \IN$, and $\lambda\in (0,\infty)$. Then, we have
\begin{eqnarray}
\hspace{1cm}		C_f(z) := |f(\omega_m(z))|^p+\lambda\bigg[B_1(f,|\omega_k(z)|) + 
A(f_0,|\omega_k(z)|)
\bigg]\leq 1 \label{eq3.3}
	\end{eqnarray}
	for $0\leq |z|= r \leq R_3$. 
	The radius $R_3$ is best possible.
\ethm	

    \br
    If we set $k=1$ and $\omega_k(z)=z$ in Theorem \ref{NEW-lem6}, then we get Theorem C. 
    \er

    The following theorem is a generalization of Theorem \ref{NEW-lem6}, which takes into account of the case $p>0$ in Theorem \ref{NEW-lem6}.

\bthm\label{NEW-lem7}
	Suppose that $f(z)=\sum_{j=0}^{\infty} a_j z^{j} \in\mathcal{B}$ with $a:=|a_0|<1$,  $\omega_n \in\mathcal{B}_n$ for $n\geq 1$, and $p > 0$, $R_4: =R^{p,m,k}_{\lambda,a}$  is
    the minimal positive root in $(0, 1)$ of the equation $\Psi_4(r)=0$, where
   \begin{equation}
   	\Psi_4(r)=(1+a r^m)^p[1+(\lambda a^2-\lambda-1)r^k]-(1-r^k)(a+r^m)^p\label{eq2.7}
   \end{equation}
for some $m, k\in \IN$ and $\lambda\in (0,\infty)$. Then   \eqref{eq3.3} holds for $0\leq |z|= r \leq R_4$. 
	The constant $R_4$ cannot be improved.
\ethm	
	
	 \br
If we set $\lambda=1$ and $a=0$ in \eqref{eq2.7}, then the minimal positive root in $(0, 1)$ of the equation $\Psi_4(r)=0$ is equivalent to determining the  minimal positive root $\beta_{k,m,p}$  in the interval $(0,1)$ of  the equation
\begin{equation}
   	   \eta_{k,m,p}(r):=1-2r^k-r^{mp}(1-r^k)=0.\label{eq2.12}
\end{equation}
We list the values of $\beta_{k,m,p}$ for certain choices of $k$, $m$ and $p$ in Table \ref{tab3}. In addition, we list below certain other simple cases.

       \begin{enumerate}

        \item  If we set $k= m=p=1$ in \eqref{eq2.12}, then it reduces to $r^2-3r+1=0$ which gives the root  $\beta_{1,1,1}=\frac{3-\sqrt{5}}{2}\approx 0.381966$.

        \item  If we allow $m\rightarrow \infty $ or $p\rightarrow \infty $ in \eqref{eq2.12}, then $\beta_{k,\infty,p}=\beta_{k,m,\infty}=\frac{1}{\sqrt[k]{2}}$. In particular, $\beta_{1,\infty,p}=\beta_{1,m,\infty}=1/2$.
        \item  If we set $\lambda=m=k=1$ with $\omega_1(z)= z$ in Theorem \ref{NEW-lem7}, we obtain \cite[Lemma 3.4]{MS2021}.

      \end{enumerate}
    \er

    \begin{table}[tbp]
\centering
\begin{tabular}{|l|l||l|l||l|l|}
\hline
$k$  &$\beta_{k,1,3}$  &$m$  &$\beta_{1,m,3}$  &$p$ &$\beta_{1,1,p}$ \\
\hline
1  &0.472213  &1  &0.472213  &1  &0.381966\\
\hline
2  &0.648791  &2  &0.496239  &2  &0.445042\\
\hline
3  &0.725563  &3  &0.499515  &3  &0.472213\\
\hline
4  &0.770224  &4  &0.499939  &4  &0.485690\\
\hline
5  &0.800095  &5  &0.499992  &5  &0.492639\\
\hline
6  &0.821776  &6  &0.499999  &6  &0.496239\\
\hline
7  &0.838383  &7  &0.500000  &7  &0.498091\\
\hline
8  &0.851600  &8  &0.500000  &8  &0.499037\\
\hline
\end{tabular}
\vspace{8pt}
\caption{Values of $\beta_{k,m,p}$ for certain choices of $k,m,p$ in Theorem \ref{NEW-lem7}} \label{tab3}
\end{table}
	
	The next theorem is a generalization of Theorem D. 
	
\bthm\label{NEW-lem8}
	Suppose that $f\in\mathcal{B}$ with $a=|f(0)|<1$ and $f(z)=\sum_{j=0}^{\infty} a_{qj} z^{qj}$ for some $q \in \IN$, $\omega_k \in\mathcal{B}_k$ for some $k\in\IN$. 
Then for $p>0$, we have the sharp inequality:
\begin{equation}
D_f(z):=a^p+\sum_{j=1}^{\infty}|a_{qj}|\,|\omega_k(z)|^{qj}+\left(\dfrac{1}{1+a}+\frac{|\omega_k(z)|^{q}}{1-|\omega_k(z)|^{q}}\right)\sum_{j=1}^{\infty}|a_{qj}|^2|\omega_k(z)|^{2qj}
\leq 1
\label{eq3.5'}
\end{equation}
for   $|z|= r \leq R_5:=R^{k,p,q}=\sqrt[qk]{\frac{\min\{p,2\}}{2+\min\{p,2\}}}$. 
\ethm	

\br
	The following observations are clear.
\begin{enumerate}
\item[{\rm (1)}]If we set $k=1$, $\omega_1(z)=z$, and $p=2$ in Theorem \ref{NEW-lem8}, then we get Theorem D(b). 

\item[{\rm (2)}] For $p=1$, by a straightforward computation, it is easy to see that \eqref{eq3.5'} holds for $r\leq \dfrac{1}{\sqrt[{qk}]{3}}$.


\end{enumerate}

\er
	
Now, we establish the following theorem as a generalization of Theorem  E 
involving another Schwarz function.

\bthm\label{NEW-lem9}
Suppose that $f \in\mathcal{B}$ has the expansion $f(z)=\sum_{j=0}^{\infty} a_j z^j$ with $a=|f(0)| < 1$, and $\omega_n \in\mathcal{B}_n$ for $n\geq 1$. Also, for $k.m\in\IN$, let
\begin{eqnarray}
E_f(z): = B_0(f,|\omega_k(z)|) + 
A(f_0,|\omega_k(z)|)
+\left|f\left(\omega_m(z)\right)-a_0\right|^2
\end{eqnarray}

Then we have the following:
\begin{enumerate}
\item[{\rm (1)}] If $k=1$, then  $E_f(z)\leq 1$
for $|z|=r \leq \frac{1}{3}$ if and only if $0\leq a\leq a^*=4\sqrt 2-5\approx0.656854$. The constant $\frac{1}{3}$ cannot be improved.

\item[{\rm (2)}]
If $2\leq k\leq \frac{m}{\log_3(2\sqrt{2}+1)}$,  and $k, m \in \IN\backslash\{1\}$,
then $E_f(z)\leq 1$ holds for $|z|=r \leq \frac{1}{\sqrt[k]3}$ if and only if $0\leq a\leq a^*=4\sqrt 2-5$, and the constant $\frac{1}{\sqrt[k]3}$ cannot be improved.
\end{enumerate}
\ethm	

   For establishing the multidimensional versions of Theorems \ref{NEW-lem4}-\ref{NEW-lem9}, we study the Bohr-type inequalities on the unit ball $B_X$ in a complex Banach space $X$. In the following, we suppose that $X$ is a complex Banach space and $\mu_n: \,B_X\rightarrow B_X$ is a Schwarz mapping such that $\mu_n$ has a zero of order $n$ at $z = 0$, where $n \in \IN$. We first state the multidimensional version of Theorem \ref{NEW-lem4}.
   For convenience, let ${\rm Hol}\,(B_X, \ID)$ denote the class of all holomorphic functions $f$ with the homogeneous polynomial expansion of the form
   \be\label{HLP-Eeq1}
   f(z)=\sum_{j=0}^\infty P_j(z),~z\in B_X
   \ee
   such that $a:=|f(0)|$.

   \bthm\label{HLP-th1}
   Let $f\in{\rm Hol}\,(B_X, \ID)$ with the form \eqref{HLP-Eeq1}. Then for $m,k\in \IN, p \in (0,2]$, we have
    \begin{align}\label{eq2.1}
   G_f(z): =&a^p+\sum_{j=1}^\infty|P_j(\mu_k(z))|+\left (\frac{1}{1+a}+ \frac{r^k}{1-r^k}\right)
    \sum_{j=1}^\infty|P_j(\mu_k(z))|^2\\
    &+\left|f\left(\mu_m(z)\right)-f(0)\right| \leq 1 \notag
   \end{align}
   for $||z||=r\leq R_1$, where $R_1:=R_{k,m}^p$ is given as in Theorem \ref{NEW-lem4}. The number $R_1$ is best possible.
   \ethm


Next, we state the multidimensional version of Theorem \ref{NEW-lem5}.

   \bthm\label{HLP-th2}
Let $f\in{\rm Hol}\,(B_X, \ID)$ with the form \eqref{HLP-Eeq1}. Then, for $m,k\in \IN, p \in (0,2]$,
$\lambda\in (0,\infty)$ and $s,\, t\in\IN$ with $0<t<s$, we have
   \begin{equation}\label{eq2.3}
   	H_f(z):=|f(\mu_m(z))|^p+\lambda \sum_{j=1}^\infty|P_{sj+t}(\mu_k(z))|\leq 1
   \end{equation}
   for $||z||=r\leq R_2$, where $R_2:=R^{p,m,k}_{\lambda,s,t}$ is as in the statement of Theorem \ref{NEW-lem5}.

   In the case when $\Psi_2(r)>0$ in some interval $(R_2, R_2+\varepsilon)$, where $\Psi_2$ is given by \eqref{eq2.4}, the number $R_2$ cannot be improved.
   \ethm

   Our next result is the multidimensional version of Theorem \ref{NEW-lem6}. 

   \bthm\label{HLP-th3}
   Let $f\in{\rm Hol}\,(B_X, \ID)$ with the form \eqref{HLP-Eeq1}.
   Then, for $m,k\in \IN $, $p \in (0,2]$ and  $\lambda\in (0,\infty)$, we have
   \begin{align}\label{eq2.5}
   	I_f(z):=|f(\mu_m(z))|^p+\lambda \left [\sum_{j=1}^\infty|P_j(\mu_k(z))|+\left(\frac{1}{1+a}+ \frac{r^k}{1-r^k}\right)
   \sum_{j=1}^\infty|P_j(\mu_k(z))|^2\right ]\leq 1
   \end{align}
   for $||z||=r\leq R_3$ where $R_3:=R_{\lambda,p}^{k,m}$ is as in the statement of Theorem \ref{NEW-lem6}. The radius $R_3$ is best possible.
    \ethm

   For our purpose, it is natural to ask about Theorem \ref{HLP-th3} for $p>2$. Finding the radius $R_3$, independent of the constant term $a=|f(0)|$, seems to be tedious. However, it is possible to state it in the following form, which is another generalization of Theorem \ref{HLP-th3}.

   \bthm\label{HLP-th4}
 Let $f\in{\rm Hol}\,(B_X, \ID)$ with the form \eqref{HLP-Eeq1}.
  Then, for $m,k\in \IN $, $p>0$ and  $\lambda\in (0,\infty)$, the inequality \eqref{eq2.5} holds for $||z||=r\leq R_{4}:=R^{k,m}_{\lambda,p,a}$, where $R_4$ is as in the statement of Theorem \ref{NEW-lem7}.
 The radius $R_{4}$ is best possible.
	\ethm

We state our next result which is a multidimensional version of  Theorem \ref{NEW-lem8}. 
	
\bthm\label{HLP-th5}
Let $f\in{\rm Hol}\,(B_X, \ID)$ with the form \eqref{HLP-Eeq1}. Then, for $m,k\in \IN $, $p>0$ and $q\in \IN$, we have
\begin{eqnarray}
J_f(z):=a^p+\sum_{j=1}^\infty|P_{qj}(\mu_k(z))| +
\left (\frac{1}{1+a}+ \frac{r^{qk}}{1-r^{qk}}\right)\sum_{j=1}^\infty|P_{qj}(\mu_k(z))|^2\leq 1\label{eq2.8}
	\end{eqnarray}
   for $||z||=r\leq R_{5}:=R^{k,p,q}$, where $R^{k,p,q}$ is as in the statement of Theorem \ref{NEW-lem8}. The radius $R_5$ is best possible.
	\ethm

Finally, we present the multidimensional version of Theorem \ref{NEW-lem9}. 

	\bthm\label{HLP-th6}
Let $f\in{\rm Hol}\,(B_X, \ID)$ with the form \eqref{HLP-Eeq1}. Then we have the following:
 \begin{enumerate}
\item[{\rm (1)}] The inequality
	 \begin{align} \label{eq2.10}
K_f(z):=&\sum_{j=0}^\infty|P_j(\mu_k(z))|+\bigg(\frac{1}{1+a}+ \frac{r^k}{1-r^k}\bigg)\sum_{j=1}^\infty|P_j(\mu_k(z))|^2\\
&+\left|f\left(\mu_m(z)\right)-f(0)\right|^2 \leq 1\notag
	 \end{align}
holds for $k=1$ with $\mu_k(z)=z$, $m\geq 1$ and $||z||=r\leq \frac{1}{3}$ if and only if $0\leq a\leq a^*=4\sqrt 2-5\approx0.656854$. The constant $\frac{1}{3}$ is best possible.
\item[{\rm (2)}] If $2\leq k\leq  \frac{m}{\log_3(2\sqrt{2}+1)}$, and $k, m \in \IN \backslash \{1\}$, then the inequality \eqref{eq2.10} holds for $||z||=r\leq\frac{1}{\sqrt[k] 3}$ if and only if $0\leq a\leq a^*=4\sqrt 2-5$.
     The constant $\frac{1}{\sqrt[k] 3}$ is best possible.
  \end{enumerate}
	\ethm

	\section{Key lemmas}\label{HLP-sec3}
	
In order to establish our main results, we need the following lemmas. 
	
	\blem\label{HLP-lem1}\cite[Proof of Theorem 1]{IRK2017}
	Suppose that $f(z)=\sum_{j=0}^{\infty} a_j z^{j} \in\mathcal{B}$ and $a:=|a_0|$. Then we have
	\begin{eqnarray*}
		\sum_{j=1}^{\infty} |a_j|r^j\leq
		\begin{cases}
			r\dfrac{1-a^2}{1-ra}\,\,\,\,\,\,\,\,\mbox{for}\,\,\,\,a\geq r,\\\\
			r\dfrac{\sqrt{1-a^2}}{\sqrt{1-r^2}} \,\,\,\,\mbox{for}\,\,\,\,a<r.
		\end{cases}
	\end{eqnarray*}
	\elem
Lemma \ref{HLP-lem1} is extracted from the two cases of the proof of Theorem 1 of \cite{IRK2017}.

	\blem\label{HLP-lem2}  \cite[Proof of Theorem 1]{SP2020}
	Suppose that $f(z)=\sum_{j=0}^{\infty} a_j z^{j} \in\mathcal{B}$. Then the following inequality holds:
	$$
B_1(f,r)+ A(f_0,r)\leq (1-|a_0|^2)\frac{r}{1-r}.
	$$
	\elem
In the above $B_1(f,r)$ and $A(f_0,r)$ are defined as in Section \ref{HLP-sec1-1}.

\blem\label{new-lem2}
Let $f\in{\rm Hol}\,(B_X, \ID)$ with the form \eqref{HLP-Eeq1}.
Then for $k \in \IN$, the following inequality
	$$
\sum_{j=1}^\infty|P_j(\mu_k(z))|+\bigg(\frac{1}{1+a}+ \frac{r^k}{1-r^k}\bigg)\sum_{j=1}^\infty|P_j(\mu_k(z))|^2\leq (1-a^2)\frac{r^k}{1-r^k}
	$$
holds for $||z||=r\in [0,1)$.
	\elem

\bpf
Let $z\in B_X \backslash \{0\}$ be fixed and $v=z/||z||$. Let
$$
F(h)=f(hv),h\in \ID.
$$
Then $F(0)=f(0)$ and $F:\,\ID \rightarrow \ID$ is holomorphic and $F(h)=\sum_{j=0}^\infty P_j(v)h^j, h\in \ID$. By Lemma \ref{HLP-lem2}, we have
\begin{eqnarray*}
\sum_{j=1}^\infty|P_j(v)|\,|h|^j+\bigg(\frac{1}{1+a}+ \frac{|h|}{1-|h|}\bigg)\sum_{j=1}^\infty|P_j(v)|^2|h|^{2j}\leq (1-a^2)\frac{|h|}{1-|h|}.
\end{eqnarray*}
Setting $h=\|z\|=r<1$, we obtain that
\begin{eqnarray*}
\sum_{j=1}^\infty|P_j(z)|+\bigg(\frac{1}{1+a}+ \frac{r}{1-r}\bigg)\sum_{j=1}^\infty|P_j(z)|^2\leq (1-a^2)\frac{\|z\|}{1-\|z\|}.
\end{eqnarray*}
In view of the inequality \eqref{eq1.2}, we obtain that
$$
\sum_{j=1}^\infty|P_j(\mu_k(z))|+\bigg(\frac{1}{1+a}+ \frac{r^k}{1-r^k}\bigg)\sum_{j=1}^\infty|P_j(\mu_k(z))|^2
\leq (1-a^2)\frac{||\mu_k(z)||}{1-||\mu_k(z)||}
\leq (1-a^2)\frac{r^k}{1-r^k}
$$
and the proof is complete.
\epf	

	\blem\label{HLP-lem3}\cite{CKX2023}
	Let $m \in \IN$ and $p\in(0,2]$. For $a\in [0,1]$, consider
	$$a\mapsto D_{p,m}(a)=\left [\left (\frac{a+r^m}{1+ar^m}\right )^p-1\right ]\varphi_0(r)+(1-a^2)N(r),
	$$
	where $\varphi_0(r)$ and $N(r)$ are some nonnegative continuous functions defined on $[0,1)$. Also, suppose that
	$$\Psi_{p,m}(r)=  p\left (\frac{1-r^m}{1+r^m}\right )\varphi_0(r) -2N(r)
	$$
	and $R:=R(m,p)$ is the minimal positive root in $(0, 1)$ of the equation $\Psi_{p,m}(r)= 0$.
	If $\Psi_{p,m}(r)\geq 0$  for $0\leq r\leq R$,  then $D_{p,m}(a)\leq 0$ for $0\leq r\leq R$.
	\elem
	
 \blem\label{new-lem3}
Let   $k, m \in \IN$ such that $k\leq \frac{m}{\log_3(2\sqrt{2}+1)}$. Also, let
\begin{eqnarray*}
\Psi_6(r) := 11r^{2m+k}-6r^{2m}-8r^{k+m}+2r^m+r^k
\end{eqnarray*}
Then  the inequality $\Psi_6(r) \geq 0$ holds for $0\leq r\leq \frac{1}{\sqrt[k]3}$.
\elem
\bpf
Note that for $r\leq \frac{1}{\sqrt[k]3}$, i.e.  $r^{-k}\geq 3$, one has $0\leq r^m<r^k\leq \frac{1}{3}$ and thus, by a straight computation, we obtain
\begin{eqnarray*}
\Psi_6(r)&=&r^k(11r^{2m}-6r^{2m-k}-8r^m+2r^{m-k}+1)\\
&=&r^k[11r^{2m}-8r^m+2r^{m-k}(1-3r^m)+1] ~\mbox{ (since $r^{-k}\geq 3$)}\\
&\geq& r^k[11r^{2m}-8r^m+6r^m(1-3r^m)+1] =r^k(-7r^{2m}-2r^m+1)\\
&=&7r^k\bigg[\frac{\sqrt {8}+1}{7}+r^m\bigg]\bigg[\frac{\sqrt {8}-1}{7}-r^m\bigg].
\end{eqnarray*}
By hypothesis, $k\leq \frac{m}{\log_3(2\sqrt{2}+1)}$ and thus, the second factor in the last relation on the right is non-negative. Indeed
$$\frac{\sqrt {8}-1}{7}-r^m\geq \frac{\sqrt {8}-1}{7}-\bigg(\frac{1}{3}\bigg)^{\frac{m}{k}}\geq  \frac{1}{2\sqrt {2}+1}-\bigg(\frac{1}{3}\bigg)^{\log_3(2\sqrt {2}+1)}= 0
$$
and the desired inequality follows.
\epf

\section{Proofs of the main results}\label{HLP-sec4}

\subsection{Proof of Theorem \ref{NEW-lem4}}	
For functions $f \in\mathcal{B}$, it has been proved by F.W. Wiener that for each integer $n\geq 1$,
we have the inequality $|a_n|\leq 1-|f(0)|^2$ (cf. \cite{IG2003}). This inequality is obtained by applying Schwarz-Pick lemma to
$$
g(z)= \frac{f(\epsilon z)+ f(\epsilon^2 z)+\cdots + f(\epsilon ^{n-1}z)+ f(z)}{n}
=  a_0+ a_nz^n+ a_{2n}z^{2n} +\cdots  ,
$$
where $\epsilon^n=1$. Again, by the classical Schwarz lemma, we have that
\be
|\omega_n(z)|\leq |z|^n
\mbox{ for $z\in \mathbb{D}$ and  for each $n\geq 1$,}
\label{eq3.1}
\ee
so that
	$$
	|f(\omega_m(z))-a_0|= \left |\sum_{j=1}^{\infty}a_j(\omega_m(z))^j\right |\leq(1-a^2)\sum_{j=1}^{\infty}r^{mj}=(1-a^2)\frac{r^m}{1-r^m}.
	$$
	Then it follows from the above inequality and Lemma \ref{HLP-lem2} that
	\begin{eqnarray*}
		A_f(z) &\leq &  a^p+(1-a^2)\frac{r^k}{1-r^k}+(1-a^2)\frac{r^m}{1-r^m} \,= \,1+(1-a^2)H(a),
	\end{eqnarray*}
	where
	$$H(a)=\frac{r^k}{1-r^k}+\frac{r^m}{1-r^m}-\frac{1-a^p}{1-a^2}.
	$$
	As $x \mapsto A(x)=\frac{1-x^p}{1-x^2}$ is decreasing on $[0,1)$ for each $p\in (0,2]$, it follows that $A(x)\geq\mathop {\lim}\limits_{a\rightarrow 1^{-}}A(x)=\frac{p}{2}$.
	Thus, $H(a)$ is obviously an increasing function of $a\in [0,1)$ and therefore,  we have
	$$
	H(a)\leq H(1)= \frac{r^k}{1-r^k}+\frac{r^m}{1-r^m}-\frac{p}{2} \leq 0,
	$$
	from which we obtain that $
	A_f(z)\leq 1$ whenever $H(1)\leq 0$, which holds for $r\leq R_1$, where $R_1$ is as in Theorem \ref{NEW-lem4}.
	
	To show that the radius $R_1$ is best possible, we consider the functions
\be\label{LSJ-eq0}
 	   \omega_n(z)=z^n ~\mbox{ for $n\geq 1$}, ~\mbox{ and }~ \varphi_a(z)=\frac{a-z}{1-az}=a-\left(1-a^{2}\right) \sum_{j=1}^{\infty}a^{j-1} z^{j}.
\ee
	Using these functions, a straightforward calculation shows that (for $z=r$)
	\begin{eqnarray*} 
		A_{\varphi_a}(z)&=& a^p+\frac{(1-a^2)r^k}{1-ar^k}+\frac{(1-a^2)^2r^{2k}}{(1+a)(1-r^k)(1-ar^k)}+(1-a^2)\frac{r^m}{1-ar^m}\\
		&=&a^p+\frac{(1-a^2)r^k}{1-r^k}+\frac{(1-a^2)r^m}{1-ar^m}\\
		&=&1+(1-a^2)\bigg(\frac{r^k}{1-r^k}+\frac{r^m}{1-ar^m}-\frac{1-a^p}{1-a^2}\bigg).
	\end{eqnarray*}
	
	Since $r\mapsto \Psi_1(r)=\frac{r^k}{1-r^k}+\frac{r^m}{1-r^m}-\frac{1-a^p}{1-a^2}$ is increasing in $(0,1)$,  $\Psi_1(r)>0$ in some interval $(R_1, R_1+\varepsilon)$, and it is easy to see that when $a\rightarrow 1^{-}$, the right hand side of the above expression is bigger than $1$.
	This verifies that the radius $R_1$ is best possible. 
\qed
\subsection{Proof of Theorem \ref{NEW-lem5}}
	Assume that $f$ satisfies the hypotheses of Theorem \ref{NEW-lem5}.
     By Schwarz-Pick lemma, we have
     \begin{eqnarray*}
        |f(u)|\leq \dfrac{|u|+a}{1+a|u|},\quad u\in \mathbb{D},
     \end{eqnarray*}
	 and, since $x\mapsto \frac{x+a}{1+ax}$ is increasing on $[0,1)$, it follows from  \eqref{eq3.1} that
	 \begin{eqnarray}
	 	|f(\omega_m(z))|\leq  \frac{|\omega_m(z)|+a}{1+a|\omega_m(z)|}\leq \frac{r^{m}+a}{1+a r^{m}},\quad |z|=r<1.
	 	\label{eq3.2}
	 \end{eqnarray}
    According to \eqref{eq3.1} and \eqref{eq3.2}, we see that $B_f(z)\leq 1+D_1(a)$, where $D_1(a)$ is as in Lemma \ref{HLP-lem3} with $ N(r)=\lambda r^{k(s+t)}/(1-r^{ks})$ and $\varphi_0(r)=1$. Now, applying  Lemma \ref{HLP-lem3}, the desired inequality  $B_f(z)\leq 1 $ holds for $0\leq r\leq R_2$, where $R_2$ is as in the statement of Theorem \ref{NEW-lem5}.
	
	To show that the  radius $R_2$ is best possible, we consider the functions
	\be\label{LSJ-eq1}
		\omega_n(z)=z^n ~\mbox{ for $n\geq 1$}, ~\mbox{ and }~ f_a(z)=\frac{z+a}{1+az}=a+\left(1-a^{2}\right) \sum_{j=1}^{\infty}(-a)^{j-1} z^{j}.
	\ee
	
	Using these functions, routine calculations show that (for $z=r$)
	\begin{eqnarray}\label{eq4.5'}
		B_{f_a}(z)=\left (\frac{a+r^m}{1+ar^m}\right )^p+\lambda(1-a^2)\frac{a^{s+t-1}r^{k(s+t)}}{1-a^sr^{ks}} =1+ (1-a)T(a),
	\end{eqnarray}
	where
	\begin{eqnarray*}
		T(a)=\frac{1}{1-a}\left[\left (\frac{a+r^m}{1+ar^m}\right )^p-1\right]+\lambda(1+a)\frac{a^{s+t-1}r^{k(s+t)}}{1-a^sr^{ks}} .
	\end{eqnarray*}
	Clearly, $B_{f_a}(z)$ is bigger than   $1$ if $T(a)>0$. In fact, for $r\in(R_2,R_2+\varepsilon)$ and $a$ close to $1$, we see that
	\begin{eqnarray*}
		\lim_{a\rightarrow 1^{-}}T(a)=2\lambda\frac{r^{k(s+t)}}{1-r^{ks}}-p\frac{1-r^m}{1+r^m}> 0.
	\end{eqnarray*}
	
	Since $2\lambda\frac{r^{k(s+t)}}{1-r^{ks}}>p\frac{1-r^m}{1+r^m}$ in some interval $(R_2,R_2+\varepsilon)$, it is easy to see that when $a\rightarrow 1^{-}$, the right hand side of the above expression in $B_{f_a}(z) > 1$. This verifies that the radius $R_2$ is best possible.
\qed

\subsection{Proof of Theorem \ref{NEW-lem6}}
As before, the hypotheses of Theorem \ref{NEW-lem6} clearly imply that
	\begin{eqnarray*}
		C_f(z)\leq \left (\frac{a+r^m}{1+ar^m}\right)^p+\lambda(1-a^2)\frac{r^k}{1-r^k}: =1 +D_2(a),
	\end{eqnarray*}
	where $D_2(a)$ is as in Lemma \ref{HLP-lem3} with $ N(r)=\lambda r^k/(1-r^k)$ and $\varphi_0(r)=1$. Now, applying  Lemma \ref{HLP-lem3}, the
	desired inequality  $C_f(z)\leq 1 $  holds for $0\leq r\leq R_3$, where $R_3$ is as in the statement of Theorem \ref{NEW-lem6}.
	
	To show that the radius $R_3$ is best possible, we consider the functions $f_a$ and $\omega_n$ given by \eqref{LSJ-eq1}.
Using these functions, straightforward calculations show that (for $z=r$)
	\begin{eqnarray*}
		C_{f_a}(z)&=&\left (\frac{a+r^m}{1+ar^m}\right )^p+\lambda \left[(1-a^2)\frac{r^k}{1-ar^k}+(1-a^2)^2\frac{r^2}{(1+a)(1-r^k)(1-ar^k)}\right]\\
		&=&\left (\frac{a+r^m}{1+ar^m}\right )^p+\lambda(1-a^2)\frac{r^k}{1-r^k}\\
		&=:&1+D_{p,m}(a).
	\end{eqnarray*}
	
	Next, we just need to show that if $r>R_3$, there exists an $a$ such that $C_{f_a}(z)$ is greater than $1$. This is equivalent to showing that $D_2(a):=D_{p,m}(a)>0$ for $r>R_3$ and for some $a\in[0,1)$.
	According to the proof of  Lemma \ref{HLP-lem3}, $a\mapsto D_2(a)$ is decreasing on $[0,1)$ for $r\in(R_3,R_3+\varepsilon)$ from which we get $D_2(a)>D_2(1)=0$.
	This verifies that the radius $R_3$ is  best possible.
\qed

\subsection{Proof of Theorem \ref{NEW-lem7}}
As in the proof of Theorem \ref{NEW-lem6}, by assumptions, \eqref{eq3.1}, \eqref{eq3.2} and Lemma \ref{HLP-lem2},
it follows easily that
\begin{eqnarray}
	   	C_f(z)&\leq& \left (\frac{a+r^m}{1+ar^m}\right)^p+\lambda(1-a^2)\frac{r^k}{1-r^k}\nonumber\\
	   	&=&1-\frac{\Psi_4(r)}{(1-r^k)(1+ar^m)^p},\label{eq3.4}
\end{eqnarray}
       where $\Psi_4(r)$ is given by \eqref{eq2.7}. Evidently \eqref{eq3.3} holds if $\Psi_4(r)\geq 0$ for $0\leq r\leq R_4$. Indeed, since $\Psi_4(0)=1-a^p>0$ and $\Psi_4(1)=\lambda(a^2-1)(1+a^p)<0$, we obtain that $\Psi_4(r)\geq 0$ holds if and only if $r\leq R_4$, where $R_4$ is the minimal positive root in $(0, 1)$ of the equation $\Psi_4(r)=0$. This gives that \eqref{eq3.3} holds for $r\leq R_4$.

       To show that the radius $R_4$ is best possible, we consider we consider the functions $f_a$ and $\omega_n$ given by \eqref{LSJ-eq1}. For these functions, direct computations yields (for $z=r$)
       	\begin{eqnarray*}
       	C_{f_a}(z)&=&\left (\frac{a+r^m}{1+ar^m}\right )^p+\lambda(1-a^2)\frac{r^k}{1-r^k}\\
        &=&1-\frac{\Psi_4(r)}{(1-r^k)(1+ar^m)^p}.
       \end{eqnarray*}
       Comparison of this expression with the right hand side of the expression in the formula \eqref{eq3.4} delivers the asserted sharpness.
\qed

\subsection{Proof of Theorem \ref{NEW-lem8}}
Clearly, the function $f$ can be represented as $f(z)=g(z^q)$, where $g\in\mathcal{B}$ and $g(z)=\sum_{n=0}^{\infty} b_nz^n$ with $b_n=a_{qn}$.
Rest of the proof follows from Lemma \ref{HLP-lem2} so that
$$
D_f(z) \leq a^p+(1-a^2)\dfrac{r^{kq}}{1-r^{kq}}= 1+\dfrac{\Psi_5(r)}{1-r^{kq}},
$$
where 
\begin{equation}
  \Psi_5(r)=(2-a^2-a^p)r^{qk}+a^p-1
\label{eq2.9}
\end{equation}
and $ \Psi_5(r)\leq 0$ is equivalent to $r^{qk}\leq \frac{1-a^p}{2-a^2-a^p}.$
By applying the result and the method of Remark 1 in \cite{SP2022}, we may verify that for $p>0$, we have
\begin{equation}
\inf\limits_{0\leq a<1}\Big\{\frac{1-a^p}{2-a^2-a^p}\Big\}=\frac{\min\{p,2\}}{2+\min\{p,2\}}.
\label{eq4.8}
\end{equation}
Hence
$$
D_f(z) \leq 1+\dfrac{\Psi_5(r)}{1-r^{kq}}\leq 1
$$
for $|z|= r \leq  R_5=R^{k,p,q}=\sqrt[qk]{\frac{\min\{p,2\}}{2+\min\{p,2\}}}$.

Moreover, the  number $R_5$ is sharp as the functions
\begin{eqnarray*}
		f^*(z)=\dfrac{a-z^q}{1-az^q}=a-(1-a^2)\sum_{j=1}^{\infty}a^{j-1}z^{qj}~\mbox{ and }~\omega_k(z)=z^k
\end{eqnarray*}
shows. In fact,
%
taking $z=r$, by a simple computation, we obtain
$$D_{f^*}(z)=a^p+(1-a^2)\frac{r^{kq}}{1-r^{kq}}.
$$
Rest of the arguments is routine and thus, $R_5$ is the best possible. \qed

\subsection{Proof of Theorem \ref{NEW-lem9}}
The case (1) follows from Theorem E 
and thus, we only need to prove the case (2).

We first consider  $r^m\leq a < 1$ and $r\leq \frac{1}{\sqrt[k]3}$. Then it follows from Lemmas \ref{HLP-lem1} and \ref{HLP-lem2} that
	\begin{align*}
		E_f(z)&\leq a+\dfrac{(1-a^2)r^k}{1-r^k}+\left(\dfrac{(1-a^2)r^m}{1-ar^m}\right)^2\nonumber
\\		&=1+\dfrac{(1-a)E_{m,k}(a,r)}{(1-r^k)(1-ar^m)^2}=:E_1(r),
	\end{align*}
	where
	\begin{eqnarray*}
		E_{m,k}(a,r)&=&(2r^k+ar^k-1)(1-ar^m)^2+r^{2m}(1-r^k)(1-a^2)(1+a)\\
		&=&-r^{2m}(1-2r^k)a^3-[r^{2m}(2-3r^{k})+2r^{m+k}]a^2\\
		&&+[r^{2m}(1-r^k)+2(1-2r^k)r^m+r^k]a+(1-r^k)r^{2m}+2r^k-1.
	\end{eqnarray*}
	Then taking the partial derivative of a yields
	\begin{eqnarray*}
		\dfrac{\partial E_{m,k}(a,r)}{\partial a}&=&-3r^{2m}(1-2r^k)a^2-2[r^{2m}(2-3r^k)+2r^{m+k}]a+r^{2m}(1-r^k)\\
&&+2(1-2r^k)r^m+r^k,\\
		\dfrac{\partial^2 E_{m,k}(a,r)}{\partial a^2}&=&-6r^{2m}(1-2r^k)a-2[r^{2m}(2-3r^k)+2r^{m+k}].
	\end{eqnarray*}
     As $\dfrac{\partial^{3} E_{m,k}(a,r)}{\partial a^{3}}\leq 0$ for $r\leq \frac{1}{\sqrt[k] 2}$,  it follows that for $a\geq r^m$
     \begin{eqnarray*}
     	\dfrac{\partial^2 E_{m,k}(a,r)}{\partial a^2}\leq \dfrac{\partial^2 E_{m,k}(a,r)}{\partial a^2}\Big|_{a=r^m}=-6r^{3m}(1-2r^k)-2r^{2m}(2-3r^k)-4r^{m+1}\leq 0, 	
     \end{eqnarray*}
     showing that $\dfrac{\partial E_{m,k}(a,r)}{\partial a}$ is a decreasing function of $a\in [r^m,1]$. Thus, we obtain
     \begin{eqnarray*}
     	\dfrac{\partial E_{m,k}(a,r)}{\partial a}\geq \dfrac{\partial E_{m,k}(a,r)}{\partial a}\Big|_{a=1}=11r^{2m+k}-6r^{2m}-8r^{k+m}+2r^m+r^k,
     \end{eqnarray*}
     which is non-negative for $r\leq \frac{1}{\sqrt[k]3}$ and  $k\leq \frac{m}{\log_3(2\sqrt{2}+1)}$ (see Lemma \ref{new-lem3}).
     Therefore, $E_{m,k}(a,r)$ is clearly monotonically increasing with respect to $a\in [r^m,1]$, and then we can obtain that
     \begin{eqnarray*}
     	E_{m,k}(a,r)\leq E_{m,k}(1,r)=(3r^k-1)(1-r^m)^2\leq 0 ~\mbox{ for }~ r\leq \frac{1}{\sqrt[k]3}.
     \end{eqnarray*}
     We conclude that $E_f(z)\leq 1$ holds for $r\leq \frac{1}{\sqrt[k]3}$ and $a\geq r^m$.

     Next, we observe that  $E_1(r)$ is an increasing function of $r\in [0,1)$ and may be written as
	\begin{eqnarray*}
		 E_1(r)&=&1+\dfrac{(1-a)}{(1-r^k)(1-ar^m)^2}\big [(2a^3+3a^2-a-1)r^{2m+k}\\&&-(a^3+2a^2-a-1)r^{2m}-(2a^2+4a)r^{m+k}+2ar^m+(a+2)r^k-1\big ].
	\end{eqnarray*}
	Now, for $r\leq \frac{1}{\sqrt[k]3}$, we have
	\begin{eqnarray*}
		E_1(r)&\leq&E_1\left (\frac{1}{\sqrt[k]3}\right )
		=1+\dfrac{(1-a)}{2(3^{\frac{m}{k}}-a)^2}\big [-a^3-3a^2+2a+2+3^{\frac{m}{k}}(1-a)(2a-3^{\frac{m}{k}})\big]\\
		&\leq&1+\dfrac{(1-a)}{2(3^{\frac{m}{k}}-a)^2} \big [-a^3-3a^2+2a+2+3(1-a)(2a-3)\big ]\\
		&=&1+\dfrac{(1-a)}{2(3^{\frac{m}{k}}-a)^2}(a^2+10a-7)\leq 1.
	\end{eqnarray*}
Here we has used the fact that  $a^2+10a-7\leq 0$ and $a\geq 0$ if and only if $0\leq a\leq a^*=4\sqrt 2-5$. In the third inequality above we have used the fact that $3^{\frac{m}{k}}(2a-3^{\frac{m}{k}})\leq 3(2a-3)$, since  $\frac{m}{k}\geq \log_3(2\sqrt{2}+1)>1$.
	
	Finally, if $0\leq a <r^m\leq (\frac{1}{3})^{\frac{m}{k}}$, then it follows from \eqref{eq3.1}, Lemmas \ref{HLP-lem2} and \ref{HLP-lem1} that
	\begin{eqnarray*}
	 	E_f(z)&\leq& a+\dfrac{(1-a^2)r^k}{1-r^k}+\dfrac{(1-a^2)r^{2m}}{1-r^{2m}}=:E_2(r).
	\end{eqnarray*}
Obviously $r\mapsto E_2(r)$ is monotonically increasing and thus, for $a <r^m\leq (\frac{1}{3})^{\frac{m}{k}}$, we see that
\begin{eqnarray*}
    	E_f(z)&\leq&E_2\left (\frac{1}{\sqrt[k] 3}\right )
        =a+\frac{1}{2}(1-a^2)+\dfrac{(\frac{1}{3})^{\frac{2m}{k}}(1-a^2)}{1-(\frac{1}{3})^{\frac{2m}{k}}}\\
    	 &\leq&\left (\frac{1}{3}\right )^{\frac{m}{k}}+\frac{1}{2}+\dfrac{(\frac{1}{3})^{\frac{2m}{k}}}{1-(\frac{1}{3})^{\frac{2m}{k}}}\\
    	&<&\frac{1}{3}+\frac{1}{2}+\frac{1}{8}<1.
    \end{eqnarray*}
    For the sharpness of $\frac{1}{\sqrt[k]3}$, we consider the functions $\omega_n$ and $\varphi_a$ given by \eqref{LSJ-eq0}.

    Taking $z=r$, by a routine computation, we obtain
    \begin{eqnarray*}
    	E_{\varphi_a}(z)
    	=a+\dfrac{(1-a^2)r^k}{1-r^k}+\dfrac{(1-a^2)^2r^{2m}}{(1-ar^m)^2}
    \end{eqnarray*}
    and the remaining arguments as before.
 \qed

\subsection{Proof of Theorem \ref{HLP-th1}}
Let $z\in B_X \backslash \{0\}$ be fixed and $v=z/||z||$. Let
$$
F(h)=f(hv),\quad h\in \ID.
$$
Then $F(0)=f(0)$ and $F:\,\ID \rightarrow \ID$ is holomorphic and $F(h)=\sum_{j=0}^\infty P_j(v)h^j, h\in \ID$.
Therefore,
$$
|F(h)-F(0)|\leq \sum_{j=1}^\infty |P_j(v)|\,|h|^j\leq (1-a^2)\frac{|h|}{1-|h|}.
$$
Setting $h=||z||=r<1$, we obtain
$$
|f(z)-f(0)|\leq \sum_{j=1}^\infty |P_j(z)|\leq (1-a^2)\frac{\|z\|}{1-\|z\|}.
$$
In view of the inequality \eqref{eq1.2}, we obtain
\begin{align}\label{eq4.1}
|f(\mu_m(z))-f(0)|\leq \sum_{j=1}^\infty |P_j(\mu_m(z))|\leq
(1-a^2)\frac{||\mu_m(z)||}{1-||\mu_m(z)||}\leq
(1-a^2)\frac{r^m}{1-r^m}.
\end{align}
By virtue of \eqref{eq4.1} and Lemma \ref{new-lem2}, we have
\begin{eqnarray}\label{eq4.2}
G_f(z)\leq a^p+(1-a^2)\frac{r^k}{1-r^k}+(1-a^2)\frac{r^m}{1-r^m}.
\end{eqnarray}
As in the proof of Theorem \ref{NEW-lem4}, it can be easily shown that the right side of the inequality \eqref{eq4.2} is not more than 1 for $r \leq R_1$, where $R_1$ is as in Theorem \ref{NEW-lem4}. Therefore, the desired inequality \eqref{eq2.1} holds for $||z||=r\leq R_1$.

Next, we prove the sharpness for $R_1$. We consider the function
$$
\eta_a(z)=\frac{a-T_v(z)}{1-aT_v(z)}=a-(1-a^2)\sum_{j=1}^\infty a^{j-1}T_v(z)^j , ~~z\in B_X, ~a \in (0,1).
$$
For $z=rv$, $v\in \partial B_X$, then $T_v(z)=r$, and
\begin{align*}
&|\eta_a(z)-\eta_a(0)|=\Big|\frac{a-T_v(z)}{1-aT_v(z)}-a\Big|=
(1-a^2)\frac{|T_v(z)|}{|1-aT_v(z)|}=(1-a^2)\frac{r}{1-ar},\\
&|P_j(z)|=(1-a^2)a^{j-1}r^j, j=1,2,...
\end{align*}
Then for $\mu_n(z)=T_v(z)^{n-1}z$, $n \in \IN$, we have
\begin{align}\label{eq4.3}
||\mu_n(z)||=|T_v(z)|^{n-1}\cdot ||z||=r^n.
\end{align}
Therefore, by a simple computation it follows that
\begin{align*}
G_{\eta_a}(z)=&a^p+(1-a^2)\sum_{j=1}^\infty a^{j-1}r^{kj}+\bigg(\frac{1}{1+a}+ \frac{r^k}{1-r^k}\bigg)(1-a^2)^2\sum_{j=1}^\infty a^{2(j-1)}r^{2kj}\\
&+(1-a^2)\frac{r^m}{1-ar^m}\\
=&1+(1-a^2)\bigg(\frac{r^k}{1-r^k}+\frac{r^m}{1-ar^m}-\frac{1-a^p}{1-a^2}\bigg).
\end{align*}
In view of the proof of Theorem \ref{NEW-lem4}, we can easily obtain that $G_{\eta_a}(z)>1$ in some interval $(R_1, R_1+\varepsilon)$. This verifies that the radius $R_1$ is best possible.
\hfill $\Box$

\subsection{Proof of Theorem \ref{HLP-th2}}
Let $z\in B_X \backslash \{0\}$ be fixed and $v=z/||z||$. Let
$$
F(h)=f(hv),h\in \ID.
$$
Then $F(0)=f(0)$ and $F:\,\ID \rightarrow \ID$ is holomorphic and $F(h)=\sum_{j=0}^\infty P_j(v)h^j, h\in \ID$.
Therefore,
\begin{eqnarray*}
\sum_{j=1}^\infty |P_{sj+t}(v)|\,|h|^{sj+t}\leq (1-a^2)\frac{|h|^{s+t}}{1-|h|^s}.
\end{eqnarray*}
By Schwarz-Pick lemma for holomorphic functions on $\ID$, we have
$$
|F(h)|\leq \frac{|F(0)|+|h|}{1+|F(0)|\,|h|}=\frac{a+|h|}{1+a|h|}.
$$
Setting $h=||z||=r<1$, we obtain
\begin{eqnarray}\label{eq4.4}
\sum_{j=1}^\infty |P_{sj+t}(z)|\leq (1-a^2)\frac{||z||^{s+t}}{1-||z||^s}
 \,\,\,\,and\,\,\,\,
|f(z)|\leq \frac{a+||z||}{1+a||z||}.
\end{eqnarray}
In view of the inequality \eqref{eq1.2} and \eqref{eq4.4}, we obtain
$$
\sum_{j=1}^\infty |P_{sj+t}(\mu_k(z))|\leq (1-a^2)\frac{||\mu_k(z)||^{s+t}}{1-||\mu_k(z)||^s}\leq (1-a^2)\frac{r^{k(s+t)}}{1-r^{ks}}
$$
and
$$
|f(\mu_m(z))|\leq \frac{a+||\mu_m(z)||}{1+a||\mu_m(z)||}\leq \frac{a+r^m}{1+ar^m}.
$$

By virtue of the two inequalities above, we have
\begin{eqnarray}\label{eq4.5}
H_f(z)\leq \bigg(\frac{a+r^m}{1+ar^m}\bigg)^p+\lambda(1-a^2)\frac{r^{k(s+t)}}{1-r^{ks}}=1+D_1(a).
\end{eqnarray}
According to the proof of Theorem \ref{NEW-lem5}, it can be easily shown that the right side of the inequality \eqref{eq4.5} is not more than 1 for $r \leq R_2$, where $R_2$ is as in Theorem \ref{NEW-lem5}. Therefore, the desired inequality \eqref{eq2.3} holds for $||z||=r\leq R_2$.

To show that the radius $R_2$ is best possible, for $z\in B_X$, $n \in \IN$ and $a \in (0,1)$, we consider the functions
\begin{align}\label{eq4.6}
\mu_n(z)=T_v(z)^{n-1}z~ \mbox{ and }~ \xi_a(z)=\frac{a+T_v(z)}{1+aT_v(z)}=a+(1-a^2)\sum_{j=1}^\infty (-a)^{j-1}T_v(z)^j.
\end{align}
Using these functions, by a routine computation it follows that
\begin{align*}
H_{\xi_a}(z)=\bigg(\frac{a+r^m}{1+ar^m}\bigg)^p+\lambda(1-a^2)\frac{a^{s+t-1}r^{k(s+t)}}{1-a^sr^{ks}}.
\end{align*}
Comparison of this expression with right side of the expression in the formula \eqref{eq4.5'} delivers the asserted sharpness.
\hfill$\Box$

\subsection{Proof of Theorem \ref{HLP-th3}}
By using the analogous proof of Theorems \ref{HLP-th1} and \ref{HLP-th2}, we may verify that
\begin{align*}
I_f(z)\leq \bigg(\frac{a+r^m}{1+ar^m}\bigg)^p+\lambda(1-a^2)\frac{r^k}{1-r^k}=1+D_2(a).
\end{align*}
In fact, according to the analogous proof of Theorem \ref{NEW-lem6}, 
we may verify that \eqref{eq2.5} holds for $r\leq R_3$, where $R_3$ is as in Theorem \ref{NEW-lem6}.

Now we prove that the constant $R_3$ is best possible. We consider the functions $\mu_n$ and $\xi_a$ given by \eqref{eq4.6}. Using these functions, straightforward calculations show that
\begin{align*}
I_{\xi_a}(z)= \bigg(\frac{a+r^m}{1+ar^m}\bigg)^p+\lambda(1-a^2)\frac{r^k}{1-r^k}=1+D_2(a).
\end{align*}
In view of the analogous proof of Theorem \ref{NEW-lem6}, we may verify that $R_3$ is the best possible. The proof is complete.
\hfill$\Box$

\subsection{Proof of Theorem \ref{HLP-th4}}
By using a similar argument of the proof of Theorem \ref{NEW-lem7} and Theorem \ref{HLP-th3}, the result can be obtained easily and we consider the same functions as in the proof of Theorem \ref{HLP-th3} to show the sharpness of $R_4$. Hence, we omit the details.
\hfill$\Box$

\subsection{Proof of Theorem \ref{HLP-th5}}
Clearly, the function $f$ can be represented as $f(z)=g(z^q)$,
where $g(z)=\sum_{j=0}^{\infty} b_jz^j$ with $b_j=a_{qj}$.
As in the proof of Theorem  \ref{HLP-th1}, we have
$$J_f(z) \leq a^p+(1-a^2)\dfrac{r^{kq}}{1-r^{kq}}= 1+\dfrac{\Psi_5(r)}{1-r^{kq}},
$$
where $\Psi_5(r)$ is given by \eqref{eq2.9}. Rest of the arguments is similar to Theorem \ref{NEW-lem8} and so skip it. Moreover, the  number $R_5$ is sharp as the functions
\begin{eqnarray*}
		\eta_a^*(z)=\dfrac{a-T_v(z)^q}{1-aT_v(z)^q}=a-(1-a^2)\sum_{j=1}^{\infty}a^{j-1}T_v(z)^{qj}~\mbox{ and }~\mu_k(z)=T_v(z)^{k-1}z
\end{eqnarray*}
shows. In fact, by a simple computation, we obtain
$$J_{\eta_a^*}(z)=a^p+(1-a^2)\frac{r^{kq}}{1-r^{kq}}.
$$
Rest of the arguments is similar and thus, $R_5$ is the best possible.
\hfill$\Box$

\subsection{Proof of Theorem \ref{HLP-th6}}
Assume the hypotheses of Theorem \ref{HLP-th6}. By means of Lemma \ref{new-lem2}, Theorem \ref{NEW-lem9} and using similar methods of the proof of Theorem \ref{HLP-th1}, we can easily obtain that \eqref{eq2.10} holds for $r\leq \frac{1}{\sqrt[k]{3}}$.

To prove the constant $\frac{1}{\sqrt[k]{3}}$ cannot be improved, in view of the analogous proof of Theorem \ref{HLP-th1}, it is sufficient to show that $\frac{1}{\sqrt[k]{3}}$ is the best possible. Here also, we consider the functions $\eta_a$ and $\mu_n$, which give the sharpness of the number $\frac{1}{\sqrt[k]{3}}$.
\hfill$\Box$

\subsection*{Acknowledgments}
This work was completed during visit of the third author to ``Lomonosov Moscow State University, Moscow Center of Fundamental and Applied Mathematics, Moscow,"  under the support of ``IMU-CDC Individual Research Travel Support".

\subsection*{Conflict of Interests}
The authors declare that they have no conflict of interest, regarding the publication of this paper.

\subsection*{Data Availability Statement}
The authors declare that this research is purely theoretical and does not associate with any data.

\end{document}